\documentclass[12pt]{amsart}
\usepackage{epsfig,amssymb, amsmath,xypic, url, amsthm, color}

\newtheorem{theorem}{Theorem}

\theoremstyle{remark}
\newtheorem{remark}{Remark}

\renewcommand{\part}{\vdash}

\renewcommand{\O}{{\mathcal O}}

\renewcommand{\part}{\vdash}

\newcommand{\PP}{\mathbf{P}}
\newcommand{\ZZ}{\mathbf{Z}}
\newcommand{\CC}{\mathbf{C}}
\newcommand{\RR}{\mathbf{R}}
\newcommand{\QQ}{\mathbf{Q}}

\begin{document}

\title{Thin monodromy in ${\rm Sp}(4)$}
\author{Christopher Brav and Hugh Thomas}

\begin{abstract}
We show that some hypergeometric monodromy groups in ${\rm Sp}(4,\ZZ)$ split as free or amalgamated products and hence by cohomological considerations give examples of  Zariski dense, non-arithmetic monodromy groups of real rank $2$. In particular, we show that the monodromy of the natural quotient of the Dwork family of quintic threefolds in $\PP^{4}$ splits as $\ZZ \ast \ZZ/5\ZZ  $. As a consequence, for a smooth quintic threefold $X$ we show that the group of autoequivalences $D^{b}(X)$ generated by the spherical twist along $\O_{X}$ and by tensoring with $\O_{X}(1)$ is an Artin group of dihedral type.
\end{abstract}

\maketitle

\section{Introduction}

The question of the arithmeticity of Zariski dense monodromy groups of families of projective varieties was raised by Griffiths-Schmid (\cite{griffithsschmid}). Non-arithmetic or `thin' examples were given by Deligne-Mostow (\cite{delignemostow})and Nori (\cite{nori}), and more general examples in a similar spirit will be provided by Fuchs-Meiri-Sarnak (\cite{fuchsmeirisarnak}). In all these examples, however, the Zariski closure of the monodromy group has real rank $1$ or is a product of groups of rank $1$, so it is natural to ask if there are examples of thin monodromy in higher rank. For more discussion of this problem, see Sarnak's Notes on thin matrix groups (\cite{sarnak}).

The simplest test cases are families over $\PP^{1} \setminus \{ {0,1, \infty} \}$ with monodromy having Zariski closure ${\rm Sp}(4,\RR)$ (which has real rank $2$). Examples of such monodromy are provided by certain hypergeometric groups  (\cite{beukersheckman}). Very recently, some of these examples have been shown to be arithmetic by Singh-Venkataramana (\cite{singhvenkataramana}), and for others
their methods are inconclusive. 

In the present paper, we focus on the examples of hypergeometric groups in ${\rm Sp}(4,\ZZ)$ having maximally unipotent monodromy at $\infty$, which have been studied by many people (\cite{chenyangyui, doranmorgan,vanvan}). In order to describe them more precisely, let us introduce some notation.

A triple $R,T,U \in {\rm GL}_{n}(\CC)$ with $R=TU$ is called {\it (irreducible) hypergeometric} if ${\rm rank}(T- I)=1$ and $R^{-1},U$ have no common eigenvalues. This second condition ensures that the given representation of the {\it hypergeometric group} $\langle R,T,U \rangle \subset {\rm GL}_{n}(\CC)$ is  irreducible. Hypergeometric groups
are precisely the monodromy groups of generalised hypergeometric ordinary differential equations (\cite{beukersheckman}). Such a group is uniquely determined, up to conjugacy, by the Jordan normal forms of $R$, $T$, and $U$, a property referred to as `rigidity'. 

We shall
be interested in the cases in which $R,T,U \in {\rm Sp}(4,\ZZ)$ and $U$ has maximally unipotent monodromy (such arise as monodromy groups of families of Calabi-Yau threefolds with $h^{2,1}=1$ over $\PP^{1} \setminus \{0,1,\infty \}$ with maximal degeneration at $\infty$). There are precisely $14$ such examples (\cite{doranmorgan}), which are labeled by a quadruple of rational numbers $(a_{1},a_{2},a_{3},a_{4})$ such that the eigenvalues of $R$ are $\exp(2 \pi i a_{j})$. In all of these cases, the group is known to be Zariski dense in ${\rm Sp}(4,\RR)$ by a criterion of Beukers-Heckman (\cite{beukersheckman}).

The $14$ possibilities are listed in the table at the end of this section. 

\begin{theorem}
In $7$ of the $14$ cases, the group splits as a free or an amalgamated product and contains a free subgroup of finite index. 
In particular, these examples give thin monodromy groups of real rank $2$.
\end{theorem}

More precisely, in the next section we show by a ping-pong argument that in $7$ cases either that there is a splitting $\langle R,T\rangle=\langle R \rangle \ast \langle T \rangle$ or $R^{k}=- I$ for some $k >0$ and there is a splitting $\langle R,T\rangle= \langle R \rangle \ast_{\langle R^{k} \rangle} \langle R^{k} ,T \rangle$. The precise form of the splitting is displayed in the table at the end of this section, from which it is immediate that either the group is free or the group generated by conjugates of $T$ by powers of $R$ is a free subgroup of finite index. 

For thinness in these cases, note that the group ${\rm Sp}(4,\ZZ)$ has cohomological dimension $2$ over $\QQ$ (\cite{leeweintraub}, Cor. 5.2.3). But since a group of finite cohomological dimension has the same cohomological dimension as its finite index subgroups (\cite{swan}, Th. 9.1), ${\rm Sp}(4,\ZZ)$ cannot have free subgroups of finite index.

Our ping-pong argument is uniform, works for $7$ of the $14$ cases, and is inconclusive for the other $7$. For $5$ of these $7$, we exhibit additional relations showing that the group does not split in the expected way. Note of course that among these $7$ other cases are $3$ examples ($(\frac{1}{6},\frac{1}{6},\frac{5}{6},\frac{5}{6})$, $(\frac{1}{6},\frac{1}{4},\frac{3}{4},\frac{5}{6})$, and $(\frac{1}{10},\frac{3}{10},\frac{7}{10},\frac{9}{10})$) of arithmetic groups from Singh-Venkataramana (\cite{singhvenkataramana}).

The monodromy groups appearing in the table are well-studied in the context of mirror symmetry and are expected to be mirror dual to certain groups of autoequivalences acting on the bounded derived category of coherent sheaves $D^{b}(X)$ on a Calabi-Yau threefold $X$. For a discussion of this, see van Enckevort-van Straten (\cite{vanvan}). 

The case with parameter $(\frac{1}{5},\frac{2}{5},\frac{3}{5},\frac{4}{5})$ corresponds to the original and most famous example of mirror symmetry, the Dwork family of quintic threefolds and its mirror family, with the autoequivalences acting on $D^{b}(X)$, $X$ a smooth quintic threefold. 
Let $\alpha=T_{\O_{X}}$ be the spherical twist functor along the structure sheaf $\O_{X}$ and $\beta=\O_{X}(1) \otimes ?$ the functor of tensoring with $\O_{X}(1)$.
We are interested in describing the subgroup of autoequivalences $\langle \alpha, \beta \rangle \subset {\rm Aut}(D^{b}(X))$. 

\begin{theorem}
The subgroup $\langle \alpha ,\beta \rangle \subset {\rm Aut}(D^{b}(X))$ is a two-generator Artin group with relation $(\alpha \beta)^{5} = (\beta \alpha)^{5}$. 
\end{theorem}

To see this, first consider the action of $\langle \alpha, \beta \rangle$ on the even dimensional cohomology $H^{{\rm ev}}(X,\QQ)$, for which it is easy to write down explicit matrices, which we denote $A$ and $B$ (\cite{vanvan}). Using this, we see that $A$ is a transvection, $B$ is maximally unipotent, and $C=AB$ has the non-trivial fifth roots of unity for eigenvalues. $A,B,C$ therefore give a hypergeometric triple and by rigidity of hypergeometric triples, we can choose bases in which $A=T$, $B=U$. 

Next, note that from the splitting $\langle A,B \rangle=\langle A, AB \rangle =\ZZ  \ast \ZZ/5\ZZ$ we see that $\langle A,B \rangle$ is generated by $A,B$ subject to the single relation $$(AB)^{5}=I.$$ 

The relation between $\alpha, \beta$ in ${\rm Aut}(D^{b}(X))$ turns out to be slightly more subtle. Namely we have $(\alpha \beta)^{5} \simeq [2]$, where $[2] \in {\rm Aut}(D^{b}(X))$ is the cohomological shift by two degrees (\cite{cankarp,kuz}). In particular, $(\alpha \beta)^{5}$ is central.  Since centrality gives $\beta (\alpha \beta)^{5} \simeq (\beta \alpha)^{5} \beta$, we have the relation

$$(\alpha \beta)^{5} \simeq (\beta \alpha)^{5}.$$
Conversely, the relation $(\alpha \beta)^{5} \simeq (\beta \alpha)^{5}$ implies centrality in a similar way.
 
It is then not hard to see that $(\alpha \beta)^{5} \simeq (\beta \alpha)^{5}$ is the only relation among $\alpha, \beta$. Indeed, consider the Artin group generated by $x,y$ subject to the relation $(xy)^{5}=(yx)^{5}$. Then we have a sequence
of surjections $\langle x, y \rangle \twoheadrightarrow \langle \alpha, \beta \rangle \twoheadrightarrow \langle A, B \rangle$, where $x \mapsto \alpha, y \mapsto \beta$ and $\alpha \mapsto A, \beta \mapsto B$. Since the only relation among the $A,B$ is $(AB)^{5}=I$, the kernel of the surjection $\langle x, y \rangle \twoheadrightarrow \langle A,B \rangle$ is generated as a normal subgroup by $(xy)^{5}$, but since $(xy)^{5}$ is central, the kernel is in fact cyclic. Similarly, the kernel of the surjection $\langle \alpha, \beta \rangle \twoheadrightarrow \langle A, B \rangle$ is cyclic and generated by $(\alpha \beta)^{5}$. The surjection $\langle x, y \rangle \twoheadrightarrow \langle \alpha, \beta \rangle$ therefore induces an isomorphism between the kernels of $\langle x, y \rangle \twoheadrightarrow \langle A, B\rangle$
and of $\langle \alpha, \beta \rangle \twoheadrightarrow \langle A, B\rangle$ and hence is itself an isomorphism.

\medskip
\medskip

We include here a table with the data for the fourteen hypergeometric groups in ${\rm Sp}(4,\ZZ)$ for which
$U$ is maximally unipotent. In the first column are the parameters $(a_{1},a_{2},a_{3},a_{4})$ giving the eigenvalues 
$\exp(2 \pi i a_{j})$ for $R$. With respect to a suitable basis (\cite{chenyangyui}), there are integers $d$ and $k$ so that we have 
\begin{equation*}
 U=\left(\begin{array}{cccc} 1&1&0&0\\
                               0&1&0&0\\
                               d&d&1&0\\
                               0&-k&-1&1 \end{array}\right),
\quad T=\left(\begin{array}{cccc} 1&0&0&0\\
                                  0&1&0&1\\
                                  0&0&1&0\\
                                  0&0&0&1 \end{array}\right) \in {\rm Sp}(4,\ZZ).
\end{equation*}    
The integers $d,k$ appear in the second two columns. In the fourth column, we describe a splitting of the group as a free or an amalgamated product in those cases in which our method is effective. In the final column we indicate cases in which
there is a relation that prevents the group from splitting in the obvious way. We do not claim in these cases to give a
complete set of relations and indeed we omit relations of the form $R^{m} =I$.

\medskip \medskip

\begin{tabular}{|c|c|c|c|c|} \hline \multicolumn{5}{|c|}{The fourteen} \\ \hline
$(a_{1},a_{2},a_{3},a_{4})$ & $d$ & $k$ & Splitting & Additional relation \\ \hline
$(\frac{1}{5},\frac{2}{5},\frac{3}{5},\frac{4}{5})$ & 5 & 5 & $\ZZ \ast \ZZ/5\ZZ$ & None \\ \hline
$(\frac{1}{8},\frac{3}{8},\frac{5}{8},\frac{7}{8})$ & 2 & 4 & $(\ZZ \times \ZZ/2\ZZ) \ast_{\ZZ/2\ZZ} \ZZ/8\ZZ $ & None \\ \hline
$(\frac{1}{12},\frac{5}{12},\frac{7}{12},\frac{11}{12})$ & 1 & 4 & $(\ZZ \times \ZZ/2\ZZ) \ast_{\ZZ/2\ZZ} \ZZ/12\ZZ$  & None \\ \hline
$(\frac{1}{2},\frac{1}{2},\frac{1}{2},\frac{1}{2})$ & 16 & 8 & $\ZZ \ast \ZZ$ & None \\ \hline
$(\frac{1}{3},\frac{1}{2},\frac{1}{2},\frac{2}{3})$ & 12 & 7 & $\ZZ \ast \ZZ$ & None \\ \hline
$(\frac{1}{4},\frac{1}{2},\frac{1}{2},\frac{3}{4})$ & 8 & 6 & $\ZZ \ast \ZZ$ & None \\ \hline
$(\frac{1}{6},\frac{1}{2},\frac{1}{2},\frac{5}{6})$ & 4 & 5 & $\ZZ \ast \ZZ$ & None \\ \hline \hline
$(\frac{1}{4},\frac{1}{3},\frac{2}{3},\frac{3}{4})$ & 2 & 3 & No & $(R^{6}T)^{2}(R^{6}T^{-1})^2$ \\ \hline 
$(\frac{1}{6},\frac{1}{6},\frac{5}{6},\frac{5}{6})$ & 1 & 2 & No & $(RT)^{8}$  \\ \hline 
$(\frac{1}{6},\frac{1}{4},\frac{3}{4},\frac{5}{6})$ & 6 & 5 & No & $(R^{6}T)^{2}(R^{6}T^{-1})^2$  \\ \hline 
$(\frac{1}{6},\frac{1}{3},\frac{2}{3},\frac{5}{6})$ & 3 & 4 & No &  $(R^{3}T)^{2}(R^{3}T^{-1})^2$  \\ \hline 
$(\frac{1}{10},\frac{3}{10},\frac{7}{10},\frac{9}{10})$  & 1 & 3 & No & $(R^{2}T)^{12}$ \\ \hline 
$(\frac{1}{4},\frac{1}{4},\frac{3}{4},\frac{3}{4})$ & 4 & 4 & ? & ? \\ \hline 
$(\frac{1}{3},\frac{1}{3},\frac{2}{3},\frac{2}{3})$ & 9 & 6 & ? & ? \\ \hline 
\end{tabular}

\medskip
\medskip

\noindent{\bf Acknowledgements}
We are grateful to Peter Sarnak, Duco van Straten, and Wadim Zudilin for their interest in this project and to Valdemar 
Tsanov for pointing out some of the relations in the above table. CB was partially supported by GRK 1463 `Analysis, Geometry, and String Theory' at Leibniz Universit\"at Hannover, MPIM Bonn, and the EPSRC Programme Grant `Motivic invariants and categorification' at Oxford University. HT was partially supported by an NSERC discovery grant and is grateful to the Hausdorff Research Institute for Mathematics and to MSRI for hospitality. We are grateful to Bielefeld University for hospitality.

\section{Ping-pong}

We begin with a quick proof via the ping-pong lemma  of the
fact that the ${\rm SL}(2)$ hypergeometric group having $R$ with parameters $(1/3,2/3)$ and $U$ unipotent splits as
$\ZZ * \ZZ/3\ZZ$.  We include details for this example because we take a very similar approach in the ${\rm Sp}(4)$ case.

By rigidity, we can choose any convenient basis in which to write $R,T,U$. Let  

$$U=\left( \begin{array}{rr} 3&4\\-1&-1 \end{array}\right),
\quad T= \left( \begin{array}{rr} 1&3\\0&1 \end{array}\right), \quad R=TU=\left( \begin{array}{rr} 0 & 1\\-1&-1 \end{array}\right).$$

Now we recall the ping-pong lemma, in a version convenient for us 
(\cite{lyndonschupp}, Prop. III.12.4):

\begin{theorem} Let a group $G$ be generated by two subgroups $G_1,G_2$, 
whose intersection is $H$.  Suppose that $G$ acts on a set $W$, and 
suppose there are disjoint non-empty subsets $X,Y$ such that 
\begin{eqnarray*}
&(G_1-H)Y \subseteq X \textrm{ and } (G_2-H)X\subseteq Y,& \\
&HY\subseteq Y \textrm{ and } HX \subseteq X.&
\end{eqnarray*}
Then $G=G_1*_{H}G_2$.  
\end{theorem}

We wish to apply the lemma with $G_1$ the subgroup generated by $T$, and
$G_2$ the subgroup generated by $R$.  We consider the given
action of $G$ on
$\RR^2$.  

Divide $\RR^2$ into twelve equally-spaced cones (with respect to the 
angle measure obtained when we put the axes at 60 degrees, as in the diagram
below):

\begin{equation*}
\begin{picture}(0,0)%
\includegraphics{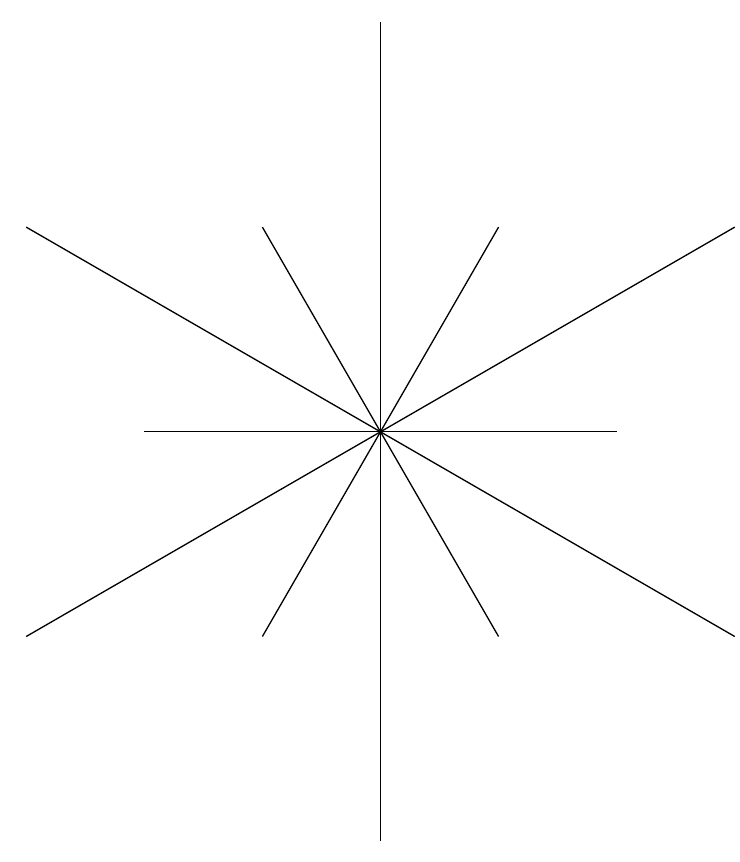}%
\end{picture}%
\setlength{\unitlength}{3315sp}%
\begingroup\makeatletter\ifx\SetFigFont\undefined%
\gdef\SetFigFont#1#2#3#4#5{%
  \reset@font\fontsize{#1}{#2pt}%
  \fontfamily{#3}\fontseries{#4}\fontshape{#5}%
  \selectfont}%
\fi\endgroup%
\begin{picture}(4212,4908)(76,-3190)
\put(2341,-3121){\makebox(0,0)[lb]{\smash{{\SetFigFont{10}{12.0}{\rmdefault}{\mddefault}{\updefault}{\color[rgb]{0,0,0}$(1,-2)$}%
}}}}
\put(1171,479){\makebox(0,0)[lb]{\smash{{\SetFigFont{10}{12.0}{\rmdefault}{\mddefault}{\updefault}{\color[rgb]{0,0,0}$(-1,1)$}%
}}}}
\put(2341,1559){\makebox(0,0)[lb]{\smash{{\SetFigFont{10}{12.0}{\rmdefault}{\mddefault}{\updefault}{\color[rgb]{0,0,0}$(-1,2)$}%
}}}}
\put(541,-691){\makebox(0,0)[lb]{\smash{{\SetFigFont{10}{12.0}{\rmdefault}{\mddefault}{\updefault}{\color[rgb]{0,0,0}$(-1,0)$}%
}}}}
\put(3421,-691){\makebox(0,0)[lb]{\smash{{\SetFigFont{10}{12.0}{\rmdefault}{\mddefault}{\updefault}{\color[rgb]{0,0,0}$(1,0)$}%
}}}}
\put( 91,-2131){\makebox(0,0)[lb]{\smash{{\SetFigFont{10}{12.0}{\rmdefault}{\mddefault}{\updefault}{\color[rgb]{0,0,0}$(-1,-1)$}%
}}}}
\put(2881,-2131){\makebox(0,0)[lb]{\smash{{\SetFigFont{10}{12.0}{\rmdefault}{\mddefault}{\updefault}{\color[rgb]{0,0,0}$(1,-1)$}%
}}}}
\put(1261,-2131){\makebox(0,0)[lb]{\smash{{\SetFigFont{10}{12.0}{\rmdefault}{\mddefault}{\updefault}{\color[rgb]{0,0,0}$(0,-1)$}%
}}}}
\put(2611,479){\makebox(0,0)[lb]{\smash{{\SetFigFont{10}{12.0}{\rmdefault}{\mddefault}{\updefault}{\color[rgb]{0,0,0}$(0,1)$}%
}}}}
\put(3871,479){\makebox(0,0)[lb]{\smash{{\SetFigFont{10}{12.0}{\rmdefault}{\mddefault}{\updefault}{\color[rgb]{0,0,0}$(1,1)$}%
}}}}
\put(181,479){\makebox(0,0)[lb]{\smash{{\SetFigFont{10}{12.0}{\rmdefault}{\mddefault}{\updefault}{\color[rgb]{0,0,0}$(-2,1)$}%
}}}}
\put(3871,-2131){\makebox(0,0)[lb]{\smash{{\SetFigFont{10}{12.0}{\rmdefault}{\mddefault}{\updefault}{\color[rgb]{0,0,0}$(2,-1)$}%
}}}}
\put(3871,-241){\makebox(0,0)[lb]{\smash{{\SetFigFont{10}{12.0}{\rmdefault}{\mddefault}{\updefault}{\color[rgb]{0,0,0}$C^-$}%
}}}}
\put(3871,-1231){\makebox(0,0)[lb]{\smash{{\SetFigFont{10}{12.0}{\rmdefault}{\mddefault}{\updefault}{\color[rgb]{0,0,0}$C^+$}%
}}}}
\put(271,-1141){\makebox(0,0)[lb]{\smash{{\SetFigFont{10}{12.0}{\rmdefault}{\mddefault}{\updefault}{\color[rgb]{0,0,0}$-C^-$}%
}}}}
\put(271,-241){\makebox(0,0)[lb]{\smash{{\SetFigFont{10}{12.0}{\rmdefault}{\mddefault}{\updefault}{\color[rgb]{0,0,0}$-C^+$}%
}}}}
\end{picture}%

\end{equation*}

We define $X$ to be the union of the interiors of the 
four cones which touch
the horizontal axis, which are labelled in the diagram as $C^+,C^-,-C^+,-C^-$.  
We define $Y$ to be the union of the interiors of
the other eight cones.  

We then must check that 
\begin{itemize}
\item [(i)] 
$X$ and $Y$ are disjoint.  
\item [(ii)] 
$R$ and $R^2$ take $X$ into $Y$.  
\item[(iii)]
$T^i$ takes $Y$ into $X$ for $i\ne 0$. 
\end{itemize}

Statements (i) and (ii) are completely trivial.  We can break statement
(iii) down into four even more trivial statements, as follows:
\begin{itemize}
\item [(iv)] $T^{-1}C^+\subseteq C^+$
\item [(v)] $TC^- \subseteq C^-$
\item [(vi)] $T^{-1}Y \subseteq \pm C^+$
\item [(vii)] $TY\subseteq \pm C^-$
\end{itemize}

This suffices.  In order to show that $T^i$ takes $Y$ into $X$ for
$i$ negative, (vi) tells us that $T^{-1}$ takes $Y$ into $\pm C^+$, and 
then (iv) tells us that applying further powers of $T^{-1}$ will not take us
outside $C^+$.  For $i$ positive, the same approach is applied using 
(vii) and (v).

\section*{Examples in ${\rm Sp}(4)$} 
We present our hypergeometric groups using the matrices

$$U=\left(\begin{array}{cccc} 1&1&0&0\\
                               0&1&0&0\\
                               d&d&1&0\\
                               0&-k&-1&1 \end{array}\right),
\quad T=\left(\begin{array}{cccc} 1&0&0&0\\
                                  0&1&0&1\\
                                  0&0&1&0\\
                                  0&0&0&1 \end{array}\right) $$
and

$$R=TU=\left(\begin{array}{cccc} 1&1&0&0\\
                                  0&1-k&-1&1\\
                                  d&d&1&0\\
                                  0&-k&-1&1 \end{array}\right),$$                                  
which preserve the standard symplectic form given by the matrix

$$J=\left(\begin{array}{cccc}0&0&1&0\\
                             0&0&0&1\\
                             -1&0&0&0\\
                             0&-1&0&0 \end{array}\right).$$

Consider the matrix:

$$B=\left(\begin{array}{cccc}-1&0&0&0\\
                             0&1&0&-1\\
                             -d&0&1&0\\
                             0&0&0&-1\end{array}\right),$$

which satisfies $B^2=I$.  We note that:

$$BRB=R^{-1} \qquad\textrm{and}\qquad BT^{-1}B=T$$

The fixed subspace of $B$ is just the span of the second and
third co-ordinate vectors.  Let us write $V$ for this subspace.  

This subspace plays a similar role in our construction to the
horizontal axis in the 2-dimensional case dicussed above, where the analogue of $B$ is a reflection 
in the horizontal axis.  

\begin{remark}
In the case that $R$ is of finite order, $V$ can be given a different 
description which is instructive, because it emphasises the analogy
of $V$ with the horizontal axis in the 2-dimensional case.  

The cases where $R$ is finite order are
exactly the cases where $R$ has four distinct eigenvalues.
In these cases, the eigenvalues are 
roots of
unity and come in complex conjugate pairs, so $R$ induces a 
decomposition of $\RR^4$ into two 2-dimensional 
subspaces $W_1$ and $W_2$, on each of which $R$ acts by a rotation in an appropriate basis.  

Each of $W_1,W_2$ has a one-dimensional 
intersection with the fixed space of $T$.  It turns
out that $V$ is the direct sum of these intersections,
and $B$ restricts to $W_i$, acting as a reflection fixing
the intersection of $W_i$ with the fixed space of $T$.  
\end{remark}

Let $$P=\log(T^{-1}R)=\log(U), \quad Q=\log(TR^{-1}).$$  

If we consider the equation $v^{T} JPv=0$, there are two 
non-zero solutions in $V$ up to scalar multiples, which are $(0,0,1,0)^{T}$ and 
$$v=(0,1,d/12-k/2,0)^{T}.$$

Now, define $C^+$ to be the open cone generated by 
$P^iv$ for $0\leq i \leq 3$ (that is to say, linear combinations of these
four vectors, with all coefficients strictly positive).  
Define $C^-=BC^+$, which can also be described as the open cone generated by
$Q^iv$ for $0 \leq i \leq 3$.  Note that $P^2v=Q^2v=(0,0,d,0)^{T}$, the
other solution to the equation which we solved for $v$.  

Finally, define $$X=\pm C^+ \cup \pm C^-, \qquad Y=\bigcup_{i \mid R^{i} \neq \pm I} R^iX.$$

We now show that $X$ and $Y$ give ping-pong tables in the first $7$ examples from the table. 

\section*{Case $(\frac{1}{5},\frac{2}{5},\frac{3}{5},\frac{4}{5})$, $d=5,k=5$}

The matrix $R$ has order 5.  $H=\{I\}$.  $G_1$ is the group generated by $T$, and 
$G_2$ is the group of order 5 generated by $R$.  

Write $M$ for the matrix whose columns
are the generating rays $v,Pv,P^{2}v,P^{3}$ of $C^+$. 
Similarly, write $N$ for the matrix whose columns are the 
generating rays $v,Qv,Q^2v,Q^3v$ of $C^-$.  We have

$$M=\left( \begin{array}{rrrr} 
     0&      1&      0&      0\\
     1&      0 &     0&      0\\
-25/12&    5/2&      5 &     0\\
     0& -25/12&      0  &   -5\end{array}\right)
\qquad N=\left( \begin{array}{rrrr}
     0&     -1&      0&      0\\
     1&  25/12&      0&      5\\
-25/12&   -5/2&      5&      0\\
     0&  25/12&      0&      5\end{array}\right)$$

The conditions we need to check are:
\begin{itemize}
\item[(i)] $X$ and $Y$ are disjoint
\item[(ii)]$(G_2-H)X\subseteq Y$
\item[(iii)] $(G_1-H)Y\subseteq X$
\end{itemize}

In order to prove (i), 
it suffices to show that $R^jC^+$ and $R^jC^-$
are disjoint from $C^+$ for $1\leq j\leq 4$. (Disjointness from
$C^-$ then follows by the symmetry with respect to $B$.)

Consider, for example,
showing that $RC^+$ is disjoint from $C^+$.  In order to do this, we expand
$M^{-1}RM$:

$$M^{-1}RM= \left(\begin{array}{rrrr}
  -23/12&   -55/12&       -5&       -5\\
       1&        1&        0&        0\\
-103/144& -131/144&   -13/12&   -25/12\\
     1/6&      1/2&        1&        1\end{array}\right)$$

Note that the entries in two of the rows are all non-negative,
while the entries in the other two rows are all non-positive. 
It follows
that the same will be true for any positive linear combination of the 
columns, and thus for any point in $RC^+$ in the basis given by the
columns of $M$.  In particular, there is
no intersection between $RC^+$ and $C^+$.

The same argument works for the other cases of (i): it is always the case
that two of the rows have non-negative entries and two of the rows have
non-positive entries.  

Condition (ii) is true by construction.  

As in the 2-dimensional case, we break condition (iii) down into four
subclaims.

\begin{itemize}
\item [(iv)] $T^{-1}C^+\subseteq C^+$.
\item [(v)] $TC^-\subseteq C^-$.  
\item [(vi)] $T^{-1}R^j(C^+\cup C^-)\subseteq \pm C^+$ for $1\leq j \leq 4$.
\item [(vii)] $TR^j(C^+\cup C^-)\subseteq \pm C^-$ for $1\leq j \leq 4$.  
\end{itemize}

Statement (iii) follows from (iv)--(vii) in essentially the same
way as in the 2-dimensional case.  
An element of $(G_1-I)Y$ is of the form 
$T^iR^jy$ for some $i\ne 0$, $1\leq j\leq 4$, $y\in Y$.  If $i<0$, then
(vi) tells us that $T^{-1}R^jy\in \pm C^+$, and then (iv) applies to tell us
that $T^iR^jy\in \pm C^+\subseteq X$.  
Similarly, if $i>0$, apply (vii) and (v).  

To establish (iv), we first argue that $T^{-1}{\overline {C^+}}\subseteq \overline{C^+}$.
To show this, it suffices to see that the generating rays of $T^{-1}{\overline{C^+}}$
are contained in $\overline{C^+}$.  We evaluate:

$$M^{-1}T^{-1}M =\left(\begin{array}{rrrr}
      1&   25/12&       0&       5\\
      0&       1&       0&       0\\
      0& 125/144&       1&   25/12\\
      0&       0&       0&       1\end{array}\right)$$

We observer that the entries are all non-negative, which is exactly what we 
needed.  
Finally, since $T^{-1}$ is an invertible linear transformation,
it takes open sets to open sets, so the image of $C^+$ will be contained in
$C^+$.  

(v) can either be shown by a completely similar argument, checking that the
entries of $N^{-1}TN$ are non-negative, 
or by deducing
it from (iv) using the symmetry encoded by $B$.

The proof of (vi) works in exactly the same way: it reduces to checking
that the entries of the following matrices are either all non-negative
or all non-positive: $M^{-1}T^{-1}R^jM$ and $M^{-1}T^{-1}R^{j}N$ 
for $1\leq j \leq 4$.  Again (vii) is established by the same argument or by
deducing it from (vi).

Ping-pong therefore establishes that the monodromy is isomorphic to
$\ZZ \ast \ZZ/5\ZZ$.

\section*{Case $(\frac{1}{8},\frac{3}{8},\frac{5}{8},\frac{7}{8})$, $d=2$, $k=4$.}
We have $R^4=-I$. Now $H=\{I,R^4\}$. 

$G_1$ is generated by $T$ and $R^4$ (forming a group isomorphic to 
$\ZZ \times \ZZ/2\ZZ$).  The subgroup $G_2$ is generated by $R$.  

Since $X$ and $Y$ are symmetric with respect to negation, the conditions
on $H$ hold.  

Condition (i) is proved exactly as before.   

As before, condition (ii), that $(G_2-H)X\subseteq Y$, is 
true by construction.  

Condition (iii) follows as before from conditions (iv)--(vii).  In conditions
(vi) and (vii), note that now $j$ runs from 1 to 3.  The conditions are 
checked the same way as before.  

Ping-pong therefore establishes that the monodromy is isomorphic to
$(\ZZ \times \ZZ/2\ZZ)*_{\ZZ/2\ZZ} \ZZ/8\ZZ$.  

\section*{Case $(\frac{1}{12},\frac{5}{12},\frac{7}{12},\frac{11}{12})$, $d=1,k=4$}

We have $R^6=-I$.  As in the previous case,
$H=\{I,-I\}$.  
$G_1$ is generated by $T$ and $R^6$ (forming a group isomorphic to 
$\ZZ \times \ZZ/2\ZZ$).  The subgroup $G_2$ is generated by $R$.

Everything is checked exactly as in the previous case.  
Ping-pong establishes that the monodromy is 
$(\ZZ \times \ZZ/2\ZZ)*_{\ZZ/2\ZZ} \ZZ/12\ZZ$.

\section*{Case $(\frac{1}{2},\frac{1}{2},\frac{1}{2},\frac{1}{2})$, $d=16,k=8$}

The order of $R$ is infinite.  
$H=\{I\}$.  $G_1$ is generated by $T$ and $G_2$ by $R$.  

The situation is very similar to the first case $d=k=5$, except
that, in order to perform the check of conditions (i), (vi), and (vii) as in 
that
case, an infinite number of checks would be required.
To establish (i) we must show that $R^jC^+$ and $R^jC^-$ are disjoint
from $C^+$ for all $j\ne 0$.  
To establish (vi), we must show
that the entries of $M^{-1} T^{-1} R^j M$ and 
$M^{-1} T^{-1} R^j N$ are all non-positive or all non-negative 
for each $j\ne 0$.  As before, (vii) then follows by symmetry.  

We note that $R$ consists of a single Jordan block with eigenvalue -1.  
We therefore define $Z= \log(-R)$, which is nilpotent.  
Now 
$$R^j=(-1)^j\exp(jZ)=(-1)^j(I +jZ +j^2Z^2/2+j^3Z^3/6).$$  
We find that $M^{-1}R^jM$ can therefore be expressed as 
$$M^{-1}R^jM= (-1)^j(I +jA_1+j^2A_2+j^3A_3)$$
for certain matrices $A_1,A_2,A_3$.  Since the entries in $A_3$ have all 
entries in two rows positive, and all entries in two rows negative, the
cubic term will eventually dominate, leading us to conclude that 
$M^{-1}R^jM$ will have all entries in two rows positive and all entries
in two rows negative, for $j$ sufficiently large.  In fact, it is 
easy to see that the cubic term dominates if $|j|\geq 6$; 
one then checks the cases 
$1\leq |j| \leq 5$ individually.  

To check condition (vi) we take a similar approach.  
We find that 
$$M^{-1}T^{-1}R^jM= (-1)^j(M^{-1}T^{-1}M +jD_1 + j^2D_2 +j^3D_3)$$
for certain matrices $D_1,D_2,D_3$.  

Since the entries of $D_3$ are all strictly negative, 
for $|j|$ sufficiently large the cubic term dominates, and the entries of 
$M^{-1}T^{-1}R^jM$ are all of the same sign.  In fact, the absolute
value of each entry of $D_3$ is at least as great as the corresponding entries
in the other matrices, so the cubic term dominates starting with $|j|\geq 6$.
We then check the cases $1\leq |j|\leq 5$ individually.  
The same analysis is applied to $M^{-1}T^{-1}R^jN$.


\section*{Case $(\frac{1}{3},\frac{1}{2},\frac{1}{2},\frac{2}{3})$, $d=12, k=7$}

As in the previous case, the order of $R$ is infinite.  
$H=\{I\}$.  $G_1$ is generated by $T$ and $G_2$ by $R$.  

The analysis is very similar to the previous case.  However, $R$ 
has eigenvalues $\exp{(\pm 2\pi i/3)}$, and
one 2-dimensional Jordan block with eigenvalue $-1$.  Therefore, we
cannot simply take the logarithm of $R$.  However, $R^6$ has all eigenvalues
equal to 1, so we can define $Z=\log(R^6)$.  Note that $Z^2=0$.  
Now 

$$R^{6n+j}=R^j\exp(nZ)=R^j(I+nZ)$$

To check that $R^mC^+$ is disjoint from $C^+$ for $m>0$, 
we observe that for each $1\leq j \leq 6$, 
$M^{-1}R^jM$ and $M^{-1}R^jZM$ have the same two rows consisting of positive
entries, and the same two rows consisting of negative entries.  This suffices,
since any positive $m$ can be written as $6n+j$ for some $1\leq j \leq 6$.  
For $m<0$,
we proceed similarly, writing $m=6n+j$ with $-6\leq j \leq -1$.  
To show that $R^kC^-$ is disjoint from $C^+$, we proceed similarly.  

To establish (vi), we consider:

$$M^{-1}T^{-1}R^{6n+j}M= M^{-1}T^{-1}R^jM + nM^{-1}T^{-1}R^j ZM$$

If $1\leq j \leq 6$, all the entries of the two matrices on the 
righthand side have the 
same sign.  This establishes that $T^{-1}R^m C^+ \subseteq C^+$ for 
any $m>0$.  

For $m<0$, we proceed similarly, writing $m=6n+j$ with $-6\leq j\leq -1$.  
Again, we check that $M^{-1}T^{-1}R^jM$ and $-M^{-1}T^{-1}R^j ZM$
have all entries with the same sign for $-6\leq j \leq -1$.  This 
establishes that $T^{-1}R^m C^+ \subseteq C^+$ for $m<0$.  

We then repeat the same two steps to show that $T^{-1}R^m C^- \subseteq C^+$.
This completes the proof of (vi).  Then (vii) follows by symmetry.  

\section*{Case $(\frac{1}{4},\frac{1}{2},\frac{1}{2},\frac{3}{4})$, $d=8,k=6$}

The matrix $R$ has eigenvalues $\pm i$ together with one Jordan block of rank two with
eigenvalue $-1$.  We
therefore proceed as in the previous case, but defining 
$Z=\log(R^4)$.  The analysis goes through in exactly the same way.

\section*{Case $(\frac{1}{6},\frac{1}{2},\frac{1}{2},\frac{5}{6})$, $d=4,k=5$}

The matrix $R$ has eigenvalues
$\exp{(\pm \pi i/3)}$ together with a Jordan block of rank two 
with eigenvalue $-1$.  We define $Z=\log(-R^3)$.  The analysis goes through
as in the previous two cases.  

\section{Logic behind the choice of ping-pong tables}
The rough outline of the shapes of the ping-pong tables was inspired from
the 2-dimensional case.  Given that, one wants to define cones 
$C^+$ and $C^-$ (as in the previous section).  
$C^+$ should be stable under $T^{-1}$ and $U$, while $C^-$ should be stable
under $T$ and $U^{-1}$.  The two cones therefore are going to lie on opposite sides
of the fixed hyperplane of $T$.  

The fact that $C^+$ and $C^-$ should have extreme rays in $V$, which is 
codimension one in the fixed hyperplane of $T$, is not obvious.  

As we already remarked, the cone $C^+$ should be 
stable under $U$.  An awkward feature
of this condition is that if we, for some reason, decide that some 
vector $x$ is in $C^+$, then it follows that the infinite set 
$Ux, U^2x,\dots$ all also lie in $C^+$.  The convex hull of this
infinite set will typically
have infinitely many extremal rays, making further analysis complicated.  

Therefore, we instead chose to require that $C^+$ be stable under $P=\log(U)$.
Since $U=\exp(P)$, stability under $P$ implies stability under $U$, but
since $P$ is nilpotent, stability under $P$ is easier to work with.  Indeed,
since $P^4=0$, for any $x$, the cone generated by $P^ix$ for $0\leq i\leq 3$
will be stable under $P$.

\bibliographystyle{plain}
\bibliography{bravthomasbiblio}

\end{document}